%&amstex          
\input amstex\documentstyle{amsppt}  
\pagewidth{12.5cm}\pageheight{19cm}\magnification\magstep1
\topmatter
\title Unipotent character sheaves and strata of a reductive group, II
\endtitle
\author G. Lusztig\endauthor
\address{Department of Mathematics, M.I.T., Cambridge, MA 02139}\endaddress    
\thanks{Supported by NSF grant DMS-2153741}\endthanks
\endtopmatter   
\document

\define\bx{\boxed}

\define\Irr{\text{\rm Irr}}

\define\si{\sim}

\define\sqc{\sqcup}

\define\op{\oplus}
   
\define\part{\partial}
\define\emp{\emptyset}

\define\m{\mapsto}
\define\do{\dots}

\define\sub{\subset}    

\define\T{\times}
\define\ti{\tilde}
\define\nl{\newline}
\redefine\i{^{-1}}

\define\un{\underline}
\define\ov{\overline}

\define\bbq{\bar{\QQ}_l}

\define\Hom{\text{\rm Hom}}
\define\End{\text{\rm End}}

\define\ind{\text{\rm ind}}

\define\g{\gamma}
\redefine\d{\delta}
\define\e{\epsilon}

\define\io{\iota}

\define\p{\pi}
\define\ph{\phi}

\define\s{\sigma}
\redefine\t{\tau}

\define\k{\kappa}

\define\kk{\bold k}

\define\CC{\bold C}

\define\NN{\bold N}

\define\QQ{\bold Q}

\define\cc{\Cal C}
\define\cd{\Cal D}

\define\cl{\Cal L}

\define\cs{\Cal S}

\define\cu{\Cal U}

\define\tH{\ti H}

\define\tL{\ti L}

\define\tW{\ti W}

\define\sha{\sharp}

\head Introduction\endhead
\subhead 0.1\endsubhead
This paper is a continuation of \cite{L22}; we extend the results of
\cite{L22} to the case of disconnected groups.
We preserve the notation of \cite{L22}.

Let $\tH$ be a reductive, possibly disconnected, group over $\CC$ with
identity component $H$ and with a given connected component $\cd$ which
generates $\tH$. Let $W$ be the Weyl group of $H$ and let
$\{s_i;i\in I\}$ be its set of simple relections;
let $\k:W@>>>W$ be the automorphism of $W$ induced by $\cd$.
Let $p$ be the order of $\k$.

Let $Pr=\{2,3,5,\do\}$ be the set of prime numbers;
let $\ov{Pr}=Pr\cup\{0\}$.
For $r\in Pr$ let $\kk_r,K_r$ be as in \cite{L22}. Let $\tH_r$
be a reductive group over $\kk_r$ with identity component $H_r$
and with a given connected component $\cd_r$ which generates $\tH_r$
such that $H_r$ has the same type as $H$ and the same
same Weyl group $W$; we also assume that the automorphism of $W$
induced by $\cd_r$ is equal to $\k$.
We set $\kk_0=K_0=\CC$, $\tH_0=\tH,H_0=H,\cd_0=\cd$.
For $r\in\ov{Pr}$ let $CS(\cd_r)$ be the (finite) set of isomorphism
classes of unipotent character sheaves on $\cd_r$. These are certain
simple perverse $K_r$-sheaves on $\cd_r$, see \cite{L09, 44.4}.
Let $Str(\cd_r)$ be the (finite) set of strata of $\cd_r$, see \cite{L21};
these are certain subsets of $\cd_r$ (unions of $H_r$-conjugacy classes of
fixed dimension) which form a partition of $\cd_r$. (These subsets are
locally closed in $\cd_r$, see \cite{C22}.)

In the remainder of this paper (except in 1.1, 3.1)
we assume that $p\in Pr$ and that either $H$ is 
quasi-simple or that $H$ is a torus.
When $H$ is not a torus then $H_{ad}$ is of type $A_n$ with $n\ge2,p=2$,
or $D_n$ with $n\ge4,p=2$, or $E_6$ with $p=2$, or $D_4$ with $p=3$;
we then say that $H$ is of type ${}^2A_n,{}^2D_n$, ${}^2E_6,{}^3D_4$
respectively.

For any $r\in\ov{Pr}$ we shall define a surjective map
$$\t:CS(\cd_r)@>>>Str(\cd_r).\tag a$$

Our definition of the map (a) is based on the generalized Springer
correspondence of \cite{L04, II} in bad characteristic.

In 1.11 we use the map (a) to give a new parametrization of
$CS(\cd_r)$. This involves associating to each stratum a 
finite group which comes from unipotent classes in bad
characteristic.

This can be also viewed as a parametrization of the unipotent
representations of the

(b) group $H_r(F_q)$ of $F_q$-rational points of a non-split form
of $H_r$ over a finite subfield $F_q$ of $\kk_r$ (with $r\in Pr$).

Here it is assumed that $H$ is not a torus and that the Frobenius map
of $H_r$ acts on $W$ as $\k$.

\subhead 0.2\endsubhead
If $G$ is
a group, $G'$ is a subgroup of $G$ and $G''$ is a subset of $G$ we set
$N_{G'}(G'')=\{g\in G';gG''g\i=G''\}$. If $g_1\in G$, then we set
$Z_{G'}(g_1)=\{g\in G';gg_1=g_1g\}$.

\head 1. The map $\t$\endhead
\subhead 1.1\endsubhead
For $r\in\ov{Pr}$ let $CS^\emp(\cd_r)$ be the subset of $CS(\cd_r)$
consisting of unipotent cuspidal character sheaves,
that is, objects $A\in CS(\cd_r)$ such that the support of $A$ is the
closure in $\cd_r$ of a single orbit of $Z^0_{H_r}\T H_r$ acting on $\cd_r$  
by $(z,g):g_1\m zgg_1g\i$; this orbit is denoted by $\s_A$. For such $A$
let $\d(A)$ be the dimension of the variety of Borel subgroups of $H_r$
that are normalized by a fixed element $h\in\s_A$ (this is independent
of the choice of $h$). We have
$$CS^\emp(\cd_r)=\sqc_{d\in\NN}CS^\emp_d(\cd_r)$$
where $CS^\emp_d(\cd_r)=\{A\in CS^\emp(\cd_r);\d(A)=d\}$.

\proclaim{Lemma 1.2} For any $d\in\NN$ the function
$r\m\sha(CS^\emp_d(\cd_r))$ from $\ov{Pr}$ to $\NN$ is constant; its
value is denoted by $N_d(\cd)\in\NN$.
\endproclaim  
This can be deduced from the explicit description of $CS(\cd_r)$ in
\cite{L09,\S46}.

\subhead 1.3\endsubhead
Let $\tW$ be the semidirect product of $W$ with a cyclic group of
order $p$ with a generator $\io$ so that in $\tW$ we have
$\io w=\k(w)\io$ for $w\in W$.

The bijection $I@>>>I$ induced by $\k:W@>>>W$ is denoted again by $\k$.
For $r\in\ov{Pr}$ and $J\sub I$ such that $\k(J)=J$ we fix a
parabolic subgroup $P_{J,r}$ and a Levi subgroup $L_{J,r}$ of $P_{J,r}$
such that $\tL_{J,r}:=N_{\tH_r}L_{J,r}\cap N_{\tH_r}P_{J,r}$
satisfies $\cd_{J,r}:=\tL_{J,r}\cap\cd_r\ne\emp$.
(For example, $L_{I,r}=H_r$ and $L_{\emp,r}$ is a maximal torus.)
We say that $J$ is $\cd$-{\it cuspidal}
if for some (or equivalently any) $r\in\ov{Pr}$ we have
$CS^\emp(\cd_{J,r})\ne\emp$. In this case $L_{J,r}$ is quasi-simple or a
torus; moreover, $J$ is uniquely determined by the type of
$(L_{J,r})_{ad}$.
Let $W_J$ be the Weyl group of $L_{J,r}$, viewed as a parabolic subgroup
of $W$.

We now fix a $\cd$-cuspidal $J$ and $A'\in CS^\emp(\cd_{J,r})$.
The induced object $\ind(A')$ is a well defined semisimple perverse
sheaf on $\cd_r$ (see \cite{L04, V, \S27}); it is in fact a direct sum of
character sheaves on $\cd_r$ (see \cite{L09, 44.10(a)}) By arguments
similar to those in \cite{L04, II, \S11}, $\End(\ind(A'))$ has a
canonical decomposition as a direct sum of lines $\op_w\cl_w$ with $w$
running through $N_{H_r}(\cd_{J,r})/L_{J,r}=N_W(\io W_J)/W_J$ such that
$\cl_w\cl_{w'}=\cl_{ww'}$ for any $w,w'$ in $N_W(\io W_J)/W_J$. (Here
$N_W(\io W_J)/W_J$ is as in 0.2 with $G=\tW,G'=W,G''=\io W_J$.) One
can verify that there is a unique $A\in CS(\cd_r)$ such that $A$ is a
summand with multiplicity one of $\ind(A')$ and the value of the
$a$-function of $W$ on the two-sided cell of $W$ attached to $A$ (see
\cite{L09, \S44}) is equal to the value of the $a$-function of $W_J$ on the
two-sided cell of $W_J$ attached to $A'$.
Now the summand $A$ of $\ind(A')$ is stable under each $\cl_w$ and we
can choose uniquely a nonzero vector $t_w\in \cl_w$ which acts on $A$
as identity. We have $t_wt_{w'}=t_{ww'}$ for any $w,w'$ in
$N_W(\io W_J)/W_J$.
We see that $\End(\ind(A'))$ is canonically the group algebra of
$N_W(\io W_J)/W_J$ (which is known to be a Weyl group). For any
$E'\in\Irr(N_W(\io W_J)/W_J)$ let $A'[E']$ be the perverse sheaf
$\Hom_{N_W(\io W_J)/W_J}(E',\ind(A'))$ on $\cd_r$. This is an object of
$CS(\cd_r)$.

\subhead 1.4\endsubhead
Let $CS'(\cd_r)$ be the set of triples $(J,E',A')$ where
$J$ is a $\cd$-cuspidal subset of $I$, $E'\in\Irr(N_W(\io W_J)/W_J)$ and
$A'\in CS^\emp(\cd_{J,r})$. We have a bijection
$$CS'(\cd_r)@>\si>>CS(\cd_r)\tag a$$
given by $(J,E',A')\m A'[E']$.

\subhead 1.5\endsubhead
Let $\cu(\cd_p)$ be the set of unipotent
$H_p$-conjugacy classes in $\cd_p$; for $\g\in\cu(\cd_p)$ the Springer correspondence
(defined in  \cite{MS04},\cite{L04, II, 11.10}),
associates to $\g$ and the constant local system $K_p$ on $\g$
an element $e(\g)\in\Irr(W^\k)$. (Here $W^\k$ is the fixed point set
of $\k:W@>>>W$; this is a Weyl group.)
Thus we have a well defined (injective) map
$e:\cu(H_p)@>>>\Irr(W^\k)$ whose image is denoted by
$\Irr_*(W^\k)$.

Let $CS^\emp(\cd_p)^{un}$ be the set of all $A\in CS^\emp(\cd_p)$ such that 
$\s_A=Z_{H_p}^0\g_A$ where $\g_A\in\cu(\cd_p)$.
Let 
$$CS'(\cd_p)^{un}=\{(J,E',A')\in CS'(\cd_p);A'\in CS^\emp(\cd_{J,p})^{un}\}.
$$

We define a map $\ti e:CS'(\cd_p)^{un}@>>>\Irr_*(W^\k)$ as follows.
Let $(J,E',A')\in CS'(\cd_p)^{un}$. Then $\g_{A'}\in\cu(\cd{J,p})$ is
defined; the restriction of $A'$ to $\g_{A'}$ is (up to a shift) a
cuspidal local system. Now the
generalized Springer correspondence in the disconnected case
  \cite{L04, II, 11.10} associates to this cuspidal local system and to
  $E'$ a unipotent $H_p$-conjugacy class $\g$ of $\cd_p$ and an
  irreducible local system on it. By definition,
  we have $\ti e(J,E',A')=e(\g)$.

  \subhead 1.6\endsubhead
Let $r\in\ov{Pr}$. In \cite{L21} a bijection
$$Str(H_r)@>>>\Irr_*(W^\k)\tag a$$
is defined. Using this and 1.4(a), we see that defining $\t$ in
0.1(a) is the same as defining a map
$$\un\t_r:CS'(\cd_r)@>>>\Irr_*(W^\k).$$

\proclaim{Lemma 1.7} Let $d\in\NN$ be such that $N_d(\cd)>0$ (see
1.2). Then one of the following holds.

(i) $CS^\emp_d(\cd_p)\sub CS^\emp(\cd_p)^{un}$. 

(ii) $CS^\emp_d(\cd_p)\cap CS^\emp(\cd_p)^{un}=\emp$. (In this case, $d=0$
and $H$ is of type ${}^2E_6$ or ${}^3D_4$.)
\endproclaim
This can be deduced from 2.6. 

\subhead 1.8\endsubhead
Let $r\in\ov{Pr}$. We will now define the map
$\un\t_r:CS'(\cd_r)@>>>\Irr_*(W^\k)$.
In the case where $I=\emp$, this map is the bijection between
two sets with one element. Assume now that $I\ne\emp$.
Let $(J,E',A')\in CS'(\cd_r)$.
We want to define $\un\t_r(J,E',A')$. 
Let $d\in\NN$ be defined by $A'\in CS^\emp_d(\cd_{J,r})$.
Assume first that we have
either $J=I$ and $d$ is as in 1.7(i) or that $J\ne I$.
 We set $\un\t_r(J,E',A')=\ti e(J,E',A')$.
 Next we assume that $J=I$ and $d$ is as in 1.7(ii). Then $d=0$
and $E'=1$. We set
$\un\t_r(I,1,A')=\text{ unit representation}$.

\subhead 1.9\endsubhead
Note that $\Irr_*(W^\k)$ can be viewed as a subset of $CS'(\cd_r)$
by $E'\m(\emp,E',K_r)$.
If $E'\in\Irr_*(W^\k)$, then from the definitions we have
$\un\t_r(\emp,E',K_r)=E'.$
It follows that $\un\t_r$ can be viewed
as a retraction of $CS'(\cd_r)$ onto its subset $\Irr_*(W^\k)$. In
particular, $\un\t_r$ is surjective.

\subhead 1.10\endsubhead
Let $\tH'_p=\tH_p/Z_{H_p}$, $H'_p=H_p/Z_{H_p}$ and let
$\p:\tH_p@>>>\tH'_p$ be the obvious map.

For any $E\in\Irr_*(W^\k)$ we denote by $\g_E$ the unique element
of $\cu_(\cd_p)$ such that $e(\g_E)=E$ (see 1.5). Let $c(E)$
be the group of connected components of the algebraic  group
$Z_{H'_p}(\p(u))$ where $u\in\g_E$; this finite group is
well defined up to isomorphism.

If $H$ is of type ${}^2A$ or ${}^2D$, then $c(E)$ is
a product of cyclic groups of order $2$.
If $H$ is of type ${}^2E_6$ then $c(E)$
is one of the following groups: $1,\cc_2,S_4$.
If $H$ is of type ${}^3D_4$ then $c(E)$
is one of $1,S_3$. (Here $\cc_m$ is a cyclic group of order $m$
and $S_n$
is the symmetric group in $n$ letters.)

We define a set $c(E)^*$ as follows.
If $H$ is of type ${}^2A$ or ${}^2D$, or if
$E\ne1$ and $H$ is of type ${}^2E_6$ or ${}^3D_4$ then
$c(E)^*$ is the set of isomorphism classes of irreducible
representations of $c(E)$ over $K_r$.
If $E=1$ and $H$ is of type ${}^2E_6$
then $c(E)^*=\sqc_{m\in\{1,2,3\}}\hat\cc_m^!$.
If $E=1$ and $H$ is of type ${}^3D_4$
then $c(E)^*=\sqc_{m\in\{1,2\}}\hat\cc_m^!$.
Here $\hat\cc_m^!$ consists of the faithful irreducible representations of
$\cc_m$.

The following theorem can be deduced from the definitions using the results in \S2.
\proclaim{Theorem 1.11}There exists a bijection
$$CS(\cd_r)@>\si>>\sqc_{E\in\Irr_*(W^\k)}c(E)^*$$
which makes the following diagram commutative:
$$\CD
CS(\cd_r)@>\si>>\sqc_{E\in\Irr_*(W^\k)}c(E)^*\\
@V\t VV           @VVV             \\
Str(\cd_r)@>\si>>\Irr_*(W^\k) \endCD$$
(The left vertical map is as in 0.1(a); the right vertical map is the
obvious one; the lower horizontal map is as in 1.6(a).)
\endproclaim

\head 2. Examples\endhead
\subhead 2.1\endsubhead
Assume that $H$ is of type ${}^2A_n$, $n\ge2$ or ${}^2D_n$, $n\ge4$.

If $H$ is of type ${}^2A_n$, let $CS''(\cd)$ be the set of pairs
$(J,E')$ where

$J$ is either $\emp$ (so that $N_W(\io W_J)/W_J=W^\k$)
or $J$ is such that $W_J$ is of type $A_{k(k+1)/2-1}$ for some 
$k\ge2$ with $k(k+1)/2-1\le n$, $k(k+1)/2-1=n\mod2$, (so that
$N_W(\io W_J)/W_J$ is a Weyl group of type $B_{(n+1-k(k+1)/2)/2}$ and
$E'\in\Irr(N_W(\io W_J)/W_J)$. (We use the convention that a Weyl
group of type $B_0$ is $\{1\}$.)

If $H$ is of type ${}^2D_n$, let $CS''(\cd)$ be the set of pairs
$(J,E')$ where 

$J$ is either $\emp$ (so that $N_W(\io W_J)/W_J=W^\k$)
 or $J$ is such that $W_J$ is of type $D_{k^2}$ for some 
odd $k\ge1$ with $k^2\le n$ (so that $N_W(\io W_J)/W_J$ is a Weyl
group of type $B_{n-k^2}$) and $E'\in\Irr(N_W(\io W_J)/W_J)$.

In any case we have a bijection $CS'(\cd_r)@>\si>>CS''(\cd)$ given by
$(J,E',A')\m(J,E')$.
Hence the map $\un\t_r$ can be viewed as a map 
$$CS''(\cd)@>>>Irr_*(W^\k).\tag a$$
Now $CS''(\cd)$ can also be viewed as the set of pairs consisting of 
a $\cd$-cuspidal $J$ and a cuspidal local system on a unipotent
$L_{J,2}$-conjugacy class in $\cd_{J,2}$.
The generalized Springer correspondence \cite{L04,II} attaches to such
a pair a unipotent $H_2$-conjugacy class in $\cd_2$ and an
irreducible local system on it.
By forgetting this last local system and by identifying $\cu(\cd_2)$
with $\Irr_*(W^\k)$ 
via $e$ (see 1.5), we obtain a map $CS''(\cd_r)@>>>\Irr_*(W^\k)$ which,
on the one hand,
is explicitly computed in \cite{L04, II} in terms of certain types of symbols and, 
on the other hand, it coincides with the map (a).

\subhead 2.2\endsubhead
In 2.3, 2.4 we describe the map $\un\t_r$ in terms of tables in the
case where $H$ is of type ${}^2E_6$ or ${}^3D_4$.
The tables are computed using results in \cite{M05}.
Our tables are almost the same as those in \cite{M05} except that
they are valid for any $r$ (not only for $r=p$ as in \cite{M05}) and
they contain some additional cuspidal objects.

In each case the table consists of a sequence
of rows. Each row represents the fibre of $\un\t_r$
over some $E\in\Irr_*(W^\k)$. The elements of that fibre are written
as symbols $(J,E',d)_{\sha=n}$. Such a symbol stands for the $n$ triples
$(J,E',A')$ in $CS'(\cd_r)$ with $J,E'$ fixed
and $A'$ running through the set
$CS^\emp_d(\cd_{J,r})$ (assumed to have $n\ge1$ elements).
When $n=1$ we omit the subscript $\sha=n$.
When $J=\emp$ we must have $d=0,n=1$ and we write $E'$ instead of
$(J,E',d)$. 
Note that the first entry in the row is $E$ itself.
Each row contains also the value of $c(E)$ (in a box).

The notation for the elements of $\Irr(W^\k)$ or
$\Irr(N_W(\io W_J)/W_J)$ is taken from \cite{M05}, but we write
$n_m$ instead of $\ph_{n,m}$.

If $H$ is of type ${}^2E_6$, we have $J=\emp$ or $W_J$ of type
$A_5$ with $N_W(\io W_J)/W_J$ of type $A_1$ or $W_J=W$ with
$N_W(\io W_J)/W_J=\{1\}$.

If $H$ is of type ${}^3D_4$, we have $J=\emp$ or $W_J=W$ with
$N_W(\io W_J)/W_J=\{1\}$.

\subhead 2.3. Table for ${}^2E_6$\endsubhead

$1_{24}$.....$\bx{1}$

$2''_{16}$.....$\bx{1}$

$4_{13}, 2'_{16}$.....$\bx{\cc_2}$

$1''_{12}$.....$\bx{1}$

$9_{10}$.....$\bx{1}$

$8'_9$.....$\bx{1}$

$8''_9$.....$\bx{1}$

$4''_7$.....$\bx{1}$

$6'_6$ .....$\bx{1}$

$9''_6,4_8$ .....$\bx{\cc_2}$

$16_5,4'_7$ .....$\bx{\cc_2}$

$12_4,6''_6,1'_{12},9'_6,(E_6,1,4)$......$\bx{S_4}$

$8'_3,(A_5,\e,0)$ .....$\bx{\cc_2}$

$8''_3,2''_4$ .....$\bx{\cc_2}$

$9_2$ .....$\bx{1}$

$4_1,2'_4$  .....$\bx{\cc_2}$

$1_0,(A_5,1,0),(E_6,1,0)_{\sha=2}$......$\bx{\cc_2}$

\subhead 2.4. Table for ${}^3D_4$\endsubhead

$1_6$ .....$\bx{1}$

$1''_3$.....$\bx{1}$

$2_2$.....$\bx{1}$

$2_1,1'_3,(D_4,1,1)$....$\bx{S_3}$

$1_0,(D_4,1,0)$....$\bx{1}$

\subhead 2.5\endsubhead
The restriction of the map $\t$ in 0.1(a) to $CS^\emp(\cd_r)$ has an
alternative definition. Namely, for $A\in CS^\emp(\cd_r)$, there is a
unique stratum $X\in Str(\cd_r)$ such that $\s_A\sub X$ (notation
of 1.1); we have $\t(A)=X$. This can be verified using the results
in 2.6.

\subhead 2.6\endsubhead
In the examples below we assume that $H$ is semisimple, $H\ne\{1\}$,
and for $A\in CS^\emp(\cd_r)$ we denote by $s$ the semisimple part of an 
element of $\s_A$. We describe the structure of $Z_{H_r}^0(s)$ in
various cases.

If $H$ is of type ${}^2A_n$ with $n=k(k+1)/2-1,k\ge2$ then:

if $r\ne 2$, $k=0\mod4$, then $Z_{H_r}^0(s)$ is of type $C_a\T D_b$
where $2a=2+4+\do+k,2b=1+3+\do+(k-1)$;  
if $r\ne 2$, $k=1\mod4$. then $Z_{H_r}^0(s)$ is of type $C_a\T B_b$
where $2a=2+4+\do+(k-1),2b+1=1+3+\do+k$;
if $r\ne 2$, $k=2\mod4$, then $Z_{H_r}^0(s)$ is of type $C_a\T B_b$
where $2a=2+4+\do+k,2b+1=1+3+\do+(k-1)$;
if $r\ne 2$, $k=3\mod4$, then $Z_{H_r}^0(s)$ is of type $C_a\T D_b$
where $2a=2+4+\do+(k-1),2b=1+3+\do+k$;
if $r=2$, then $Z_{H_r}^0(s)=H_r$.

If $H$ is of type ${}^2D_{k^2}$ with $k\in\{3,5,7,\do\}$ then:

if $r\ne 2$, then $Z_{H_r}^0(s)$ is of type $B_a\T B_a$ with $2a+1=k^2$;
if $r=2$, then $Z_{H_r}^0(s)=H_r$.

If $H$ is of type ${}^2E_6$, then:

if $d=4$ and $r\ne2$, then $Z_{H_r}^0(s)$ is of type $F_4$;
if $d=4$ and $r=2$, then $Z_{H_r}^0(s)=H_r$;
if $d=0$ and $r\ne3$, then $Z_{H_r}^0(s)$ is of type $A_2\T A_2\T A_2$;
if $d=0$ and $r=3$, then $Z_{H_r}^0(s)$ is of type $F_4$.

If $H$ is of type ${}^3D_4$, then:

if $d=1$ and $r\ne3$, then $Z_{H_r}^0(s)$ is of type $G_2$;
if $d=1$ and $r=3$, then $Z_{H_r}^0(s)=H_r$;
if $d=0$ and $r\ne2$, then $Z_{H_r}^0(s)$ is of type
$A_1\T A_1\T A_1\T A_1$;
if $d=0$ and $r=2$, then $Z_{H_r}^0(s)$ is of type $G_2$.

\head 3. Concluding remarks\endhead
\subhead 3.1\endsubhead
In this section we assume that $r\in\ov{Pr}$ and that $\tH=H=\cd$.
For any stratum $S$ of $H_r$ we denote by $[S]$ the union of strata $S'$
of $H_r$ such that $S\le S'$ for the natural partial order on $Str(H_r)$.
(In this partial order the set of regular elements form the biggest stratum
and the center forms the smallest stratum.)
We would like to give an alternative (conjectural) description of
our map $\t:CH(H_r)@>>>Str(H_r)$. Let $A\in CH(H_r)$.
Ler $\cs(A)$ be the set of strata $S$ of $H_r$ such that for some integer $n$,
the restriction $A[n]|_{[S]}$ is constructible. We conjecture that

(a) $\cs(A)$ has a unique minimal element $S_A$ (for $\le$) and that
$S_A=\t(A)$.
\nl
This is correct if $A$ is cuspidal.
If $A$ is the constant sheaf $\bbq$ (up to shift) then $\cs(A)=Str(H_r)$ and
$S_A$ is the centre of $H_r$. If $A$ is the unipotent character sheaf on $H_r$
corresponding to the sign representation of $W$ then $\cs(A)$ consists of
a single element (the stratum of regular elements) and (a) is again correct.

A similar conjecture can be given 
without the assumption $\tH=H$.

\enddocument